\documentclass[11pt]{amsart}
\usepackage{a4}
\usepackage{latexsym}
\usepackage{amsfonts}
\usepackage{amsmath}
\usepackage{amssymb}
\usepackage{amscd}
\usepackage{amsthm}
\theoremstyle{plain}

   \newtheorem{theorem}                    {Theorem}
    \newtheorem{lemma}      [theorem]       {Lemma}
    \newtheorem{corollary}  [theorem]       {Corollary}
    \newtheorem{proposition}[theorem]       {Proposition}

\begin{document}

\catcode`@=11
\atdef@ I#1I#2I{\CD@check{I..I..I}{\llap{$\m@th\vcenter{\hbox
  {$\scriptstyle#1$}}$}
  \rlap{$\m@th\vcenter{\hbox{$\scriptstyle#2$}}$}&&}}
\atdef@ E#1E#2E{\ampersand@
  \ifCD@ \global\bigaw@\minCDarrowwidth \else \global\bigaw@\minaw@ \fi
  \setboxz@h{$\m@th\scriptstyle\;\;{#1}\;$}%
  \ifdim\wdz@>\bigaw@ \global\bigaw@\wdz@ \fi
  \@ifnotempty{#2}{\setbox@ne\hbox{$\m@th\scriptstyle\;\;{#2}\;$}%
    \ifdim\wd@ne>\bigaw@ \global\bigaw@\wd@ne \fi}%
  \ifCD@\enskip\fi
    \mathrel{\mathop{\hbox to\bigaw@{}}%
      \limits^{#1}\@ifnotempty{#2}{_{#2}}}%
  \ifCD@\enskip\fi \ampersand@}
\catcode`@=\active

\renewcommand{\labelenumi}{\alph{enumi})}
\newcommand{\chr}{\operatorname{char}}
\newcommand{\isom}{\stackrel{\sim}{\longrightarrow}}
\newcommand{\Aut}{\operatorname{Aut}}
\newcommand{\Hom}{\operatorname{Hom}}
\newcommand{\End}{\operatorname{End}}
\newcommand{\HOM}{\operatorname{{\mathcal H{\mathfrak{om}}}}}
\newcommand{\EXT}{\operatorname{\mathcal E{\mathfrak xt}}}
\newcommand{\Tot}{\operatorname{Tot}}
\newcommand{\Ext}{\operatorname{Ext}}
\newcommand{\Gal}{\operatorname{Gal}}
\newcommand{\Pic}{\operatorname{Pic}}
\newcommand{\Spec}{\operatorname{Spec}}
\newcommand{\trdeg}{\operatorname{trdeg}}
\newcommand{\im}{\operatorname{im}}
\newcommand{\coim}{\operatorname{coim}}
\newcommand{\coker}{\operatorname{coker}}
\newcommand{\gr}{\operatorname{gr}}
\newcommand{\id}{\operatorname{id}}
\newcommand{\Br}{\operatorname{Br}}
\newcommand{\cd}{\operatorname{cd}}
\newcommand{\CH}{CH}
\newcommand{\Alb}{\operatorname{Alb}}
\renewcommand{\lim}{\operatornamewithlimits{lim}}
\newcommand{\colim}{\operatornamewithlimits{colim}}
\newcommand{\rk}{\operatorname{rank}}
\newcommand{\codim}{\operatorname{codim}}
\newcommand{\NS}{\operatorname{NS}}
\newcommand{\cone}{{\rm cone}}
\newcommand{\rank}{\operatorname{rank}}
\newcommand{\ord}{{\rm ord}}
\newcommand{\f}{{\mathcal G}}
\newcommand{\g}{{\mathcal F}}
\renewcommand{\L}{{\mathcal L}}
\newcommand{\du}{{\mathcal D}}
\newcommand{\G}{{\mathbb G}}
\newcommand{\N}{{\mathbb N}}
\newcommand{\A}{{\mathbb A}}
\newcommand{\Z}{{{\mathbb Z}}}
\newcommand{\Q}{{{\mathbb Q}}}
\newcommand{\R}{{{\mathbb R}}}
\newcommand{\B}{{\mathbb Z}^c}
\renewcommand{\H}{{{\mathbb H}}}
\renewcommand{\P}{{{\mathbb P}}}
\newcommand{\F}{{{\mathbb F}}}
\newcommand{\m}{{\mathfrak m}}
\newcommand{\Sch}{{\text{\rm Sch}}}
\newcommand{\et}{{\text{\rm et}}}
\newcommand{\fl}{{\text{\rm fl}}}
\newcommand{\eh}{{\text{\rm eh}}}
\newcommand{\Zar}{{\text{\rm Zar}}}
\newcommand{\Nis}{{\text{\rm Nis}}}
\newcommand{\tr}{\operatorname{tr}}
\newcommand{\tor}{{\text{\rm tor}}}
\newcommand{\red}{{\text{\rm red}}}
\newcommand{\sn}{{+}}
\newcommand{\Div}{\operatorname{Div}}
\newcommand{\Ab}{{\text{\rm Ab}}}
\newcommand{\DD}{{\mathbb Z}^c}
\renewcommand{\div}{\operatorname{div}}
\newcommand{\corank}{\operatorname{corank}}
\renewcommand{\O}{{\mathcal O}}
\newcommand{\C}{{\mathbb C}}
\newcommand{\p}{{\mathfrak p}}
\renewcommand{\proof}{\noindent{\it Proof. }}
\newcommand{\proofend}{\hfill $\Box$ \\}
\newcommand{\rem}{\noindent {\it Remark. }}
\newcommand{\example}{\noindent {\bf Example. }}
\newcommand{\ar}{{\text{\rm ar}}}
\newcommand{\del}{{\delta}}

\interfootnotelinepenalty=10000

\title{\large The affine part of the Picard scheme (corrected).}
\author{\large Thomas Geisser}
\address{Dep.\ of Math.\\
Rikkyo University\\
Japan}

\thanks{Supported in part by NSF grant No.0556263}

\begin{abstract}\normalsize
We describe the maximal torus and maximal unipotent subgroup of
the Picard variety of a proper scheme over a perfect field.
(This is a corrected and improved version of the article originally published
in Comp. Math. 145 (2009)).
\end{abstract}

\maketitle

\section{Introduction}
For a proper scheme $p:X\to k$ over a perfect field, the Picard scheme
$\Pic_X$ representing the functor  $T\mapsto H^0(T_\et,R^1p_*\G_m)$
exists, and its connected component $\Pic^0_X$
is separated and of finite type \cite[II 15]{murre}. 
By Chevalley's structure theorem
\cite{chevalley},
the reduced connected component $\Pic_X^{0,\red}$
is an extension of an abelian variety $A_X$ by a linear algebraic group $L_X$:
\begin{equation}\label{chev}
0\to L_X \to \Pic^{0,\red}_X \to A_X \to 0.
\end{equation}
The commutative, smooth affine group scheme $L_X$ is the
direct product of a torus $T_X$ and a unipotent group $U_X$.
The following theorem completely characterizes $T_X$:

\begin{theorem}\label{torus}
If $X$ is proper over a perfect field, then the cocharactermodule
$\Hom_{\bar k}(\G_m,T_X)$ of the maximal torus of $\Pic_X$ is
isomorphic to $H^1_\et(\bar X,\Z)$ as a Galois-module.
\end{theorem}

To analyze the unipotent part, we let $\Pic(X[t])_{[1]}$ be the
typical part, i.e. the
subgroup of elements $x$ of $\Pic(X[t])$ such that the map $X[t]\to
X[t]$, $t\mapsto nt$ sends $x$ to $nx$.

\begin{theorem}\label{npic}
Let $X$ be proper over a perfect field. Then
$\Pic(X[t])_{[1]}$ is isomorphic to the group of
morphisms of schemes $f:\G_a\to U_X$ satisfying $f(nx)=nf(x)$ for every
$n\in \Z$.
In particular, $\Hom_k(\G_a,U_X)\subseteq \Pic(X[t])_{[1]}$, and this
is an equality in characteristic $0$.
\end{theorem}

To get another description of $U_X$, we assume that $X$ is reduced
(the map on the Picard scheme induced by the map
$X^\red\to X$ is well understood by the work of Oort \cite{oortthesis}).
The semi-normalization $X^\sn\to X$ is
the largest scheme between $X$ and its normalization which is strongly
universally homeomorphic to $X$ in the sense that the map $X^\sn\to X$
induces an isomorphism on all residue fields. A Theorem of Traverso \cite{traverso}
implies that $\Pic(X[t])_{[1]}$, hence $U_X$, vanishes if $X$ is reduced
and seminormal. We use this to show

\begin{theorem}\label{unipotent}
Let $X$ be reduced and proper over a perfect field.

a) We have a short exact sequence 
\begin{equation}\label{ssss}
0\to K_X\to \Pic_X^{0,\red} \to \Pic_{X^\sn}^{0,\red}\to 0,
\end{equation}
and inclusions of unipotent group schemes 
$$ U_X\subseteq K_X \subseteq p_*(\G_{m,X^\sn}/\G_{m,X})$$ 
with quotients finite $p$-primary group schemes.

b)  The group scheme $p_*(\G_{m,X^\sn}/\G_{m,X})$ represents the functor
$$ T\mapsto \{\O_{X\times T}\text{-line bundles }
{\mathcal L}\subseteq \O_{X^\sn\times T}\text{ which are invertible in }
\O_{X^\sn\times T} \}.$$
\end{theorem}

{\it Notation:} For a field $k$, we denote by $\bar k$ its
algebraic closure, and for a scheme $X$ over $k$ we let
$\bar X=X\times_k\bar k$. Unless specified otherwise, all
extension and homomorphism groups are considered on the fpqc site.
\footnote{In \cite{ich} we used the \'etale topology}

\medskip

{\sc Acknowledgements: } 

This (original) paper was written while the author was
visiting T.\ Saito at the University of Tokyo, whom we thank for his
hospitality. We are indebted to G.\ Faltings for pointing out a mistake
in a previous version, and the referee, 
whose comments helped to improve the exposition and to give more 
concise proofs.
O.\ Gabber pointed out mistakes in the original version and suggested improvements.

\section{The torus}

\begin{proposition}\label{bef}
If $p:X\to k$ is reduced, geometrically connected, and proper over a perfect field,
then $\G_{m,k}\to p_*\G_{m,X}$ is an isomorphism. Moreover, if 
$f:X'\to X$ is a universal homeomorphism and $X'$ is reduced as well, 
then $f$ induces an isomorphism $p_*\G_{m,X}\cong p_*\G_{m,X'}$.
\footnote{This replaces \cite[Prop.\ 9 a)]{ich} which is incorrect as stated
because the induction step in the proof does not preserve the hypothesis on reducedness.}
\end{proposition}

\proof
Since any scheme $T$ over $k$ is flat, we have by flat base change
$R^jq_* \O_{X_T}= H^j(X,\O_X)\otimes_k \O_T$, where $q:X_T\to T$
is the projection. In particular,
$$p_*\G_{m,X}(T):=\Gamma(X\times T,\O_{X\times T})^\times
=(\Gamma(X,\O_X)\otimes \Gamma(T,\O_T))^\times,$$
and it suffices to show that $\Gamma(X,\O_X)\cong \Gamma(X',\O_{X'})\cong k$.
Since $\Gamma(\bar X,\O_{\bar X})^{Gal(\bar k/k)}=\Gamma(X,\O_X)$,
we can assume that $k$ is algebraically closed and that $X$ is connected,
in which case the statement follows
because $X$ and $X'$
are reduced, proper, connected, and have a $k$-rational point.
\proofend

\begin{lemma}
For any scheme $X$ we have isomorphisms 
$$H^1_\et(X,\Z)\cong H^1_\fl(X,\Z)\cong 
\Ext^1_X(\G_{m,X},\G_{m,X}).$$
\end{lemma}

\proof 
The first isomorphism is \cite[III Rem.\ 3.11(b)]{milnebook}.
To prove the second isomorphism,  we note that
$\HOM_X(\G_{m,X},\G_{m,X})\cong \Z_X$ by \cite[VIII Cor.\ 1.5]{sga3},
and that $\EXT^1_X(\G_{m,X},\G_{m,X})$ is isomorphic to the group of 
extensions of group schemes \cite[Cor.\ 17.5]{oort}, which vanishes by 
\cite[VIII Prop.\ 3.3.1]{sga7}.\footnote{This was claimed without
proof in \cite{ich}.}
Hence we obtain the isomorphism from the spectral sequence 
\cite[III Thm.1.22]{milnebook}
$$ E_2^{s,t}=H^s_\fl(X,\EXT^t_X(\G_{m,X},\G_{mX}))\Rightarrow 
\Ext^{s+t}_X(\G_{m,X},\G_{m,X}).$$
\proofend

\proof (Theorem \ref{torus})
Since the maps defined below are natural, we can assume that $k$ is 
algebraically closed and $X$ is connected.
We can also assume that $X$ is reduced, because
$H^1_\et(X,\Z)\xrightarrow{\sim} H^1_\et(X^\red,\Z)$,
and the map $\Pic_X \to \Pic_{X^\red}$ has unipotent kernel
and cokernel \cite[Cor.\ page 9]{oortthesis}.
It suffices to calculate $\Hom_k(\G_{m,k},\Pic_X)$,
because there are no homomorphisms from $\G_m$ to commutative group 
schemes other than tori \cite[p. 81]{oort}. By Yoneda's Lemma, the latter group
is isomorphic to the group of homomorphisms of sheaves on the fpqc
site $\Hom_k(\G_{m,k},R^1p_*\G_{m,X})$. The Leray spectral sequence
\begin{equation}\label{sse}
E_2^{s,t}=\Ext^s_k(\G_{m,k},R^tp_*\G_{m,X})\Rightarrow \Ext^{s+t}_X(\G_{m,X},\G_{m,X}).
\end{equation}
gives an exact sequence
$$ 0\to \Ext^1_k(\G_{m,k},p_*\G_{mX})\to
\Ext^1_X(\G_{m,X},\G_{m,X})\to 
\Hom_k(\G_{m,k},R^1p_*\G_{m,X})\stackrel{\delta_X}{\longrightarrow}
\Ext^2_k(\G_{m,k},p_*\G_{m,X}).$$
By Proposition \ref{bef} the left term agrees with $\Ext^1_k(\G_{m,k},\G_{m,k})$,
and this vanishes by \cite[Cor.\ 17.5]{oort}.
Thus it suffices to show that $\delta_X$ is the zero map
\footnote{The remainder of the proof is a simplification suggested
by O.\ Gabber.}.
Choose a closed point of $X$ and
let $i:Z\to X$ be the corresponding closed subscheme. Since 
$p_*\circ i_*=\id$ we have $R^sp_*i_*=R^s(p\circ i)_* =0$
for $s>0$. Hence we obtain a diagram
$$\begin{CD}
\Hom_k(\G_{m,k},R^1p_*\G_{m,X})@>>> \Hom_k(\G_{m,k},R^1p_*i_*\G_{m,Z})=0\\
@V\delta_XVV @V\delta_ZVV \\
\Ext^2_k(\G_{m,k},p_*\G_{m,X})@>\sim>> \Ext^2_k(\G_{m,k},\G_{m}).
\end{CD}$$
By Proposition \ref{bef}, the lower horizontal map is an isomorphism.
\proofend

{\it Remark.}
The example in \cite[Prop. 8.2]{ichweilII} shows that the map
$H^i_\et(\bar X,\Z)\to \Ext^i_{\bar X}(\G_m,\G_m)$ is not an
isomorphism for $i\geq 2$. One can ask if it is an
isomorphism if one replaces $H^i_\et(\bar X,\Z)$ by the eh-cohomology
group $H^i_{eh}(\bar X,\Z)$ of \cite{ichweilII}.

\smallskip

\noindent {\it Example.}
If $X$ is the node over an algebraically
closed field, then $H^1_\et(X,\Z)\cong \Z$, and $T_X\cong \G_m$.
Let $X$ be a node with non-rational tangent slopes at the singular point.
Base changing to the
algebraic closure, one sees that $H^1_\et(\bar X,\Z)\cong \Z$,
with Galois group acting as multiplication by $-1$, hence
$T_X$ is an anisotropic torus.


Using the theorem, we are able to recover the torsion of
$T_X$, $A_X$ and the diagonalizable part of $NS_X:=\Pic_X/\Pic_X^{0,\red}$
in terms of etale cohomology:

\begin{corollary}\label{detail}
Let $X$ be proper over a perfect field $k$. Then we have canonical isomorphisms
\begin{align*}
H^1_\et(\bar X,\Z)\otimes\Q/\Z &\cong \colim \Hom_{\bar k}(\mu_m,T_X);\\
\Div({}_\tor H^2_\et(\bar X,\Z)) &\cong \colim \Hom_{\bar k}(\mu_m,A_X);\\
{}_\tor H^2_\et(\bar X,\Z)/\Div &\cong \colim \Hom_{\bar k}(\mu_m,NS_X).
\end{align*}
\end{corollary}

\proof
Taking the colimit of the isomorphism
$H^1_\et(\bar X,\Z/m)\cong \Hom_{\bar k}(\mu_m,\Pic_X)$
of \cite[Prop.4.16]{milnebook} or \cite[\S 6.2]{raynaud},
we obtain
$ H^1_\et(\bar X,\Q/\Z)\cong \colim \Hom_{\bar k}(\mu_m,\Pic_X)$. 
Since $\Ext^1_{\bar k}(\G_m,T_X)=0$, Theorem \ref{torus} implies
that $\Hom_{\bar k}(\mu_m,T_X)\cong \Hom_{\bar k}(\G_m,T_X)/m\cong
H^1_\et(\bar X,\Z)/m$.
Consider the commutative diagram:
$$\begin{CD}
\colim \Hom_{\bar k}(\mu_m,T_X)@= H^1_\et(\bar X,\Z)\otimes \Q/\Z\\
@VVV @VVV \\
\colim \Hom_{\bar k}(\mu_m,\Pic^{0,\red}_X)@>>>
\colim \Hom_{\bar k}(\mu_m,\Pic_X)
@>>> \colim \Hom_{\bar k}(\mu_m,NS_X)\\
@VVV @VVV @|\\
\colim \Hom_{\bar k}(\mu_m,A_X)@>f>> {}_\tor H^2_\et(\bar X,\Z)
@>>> \coker f.
\end{CD}$$
The middle column is the short exact coefficient sequence.
The left column and middle row are short exact because $\Ext^1_{\bar k}(\mu_m,T_X)
=\Ext^1_{\bar k}(\mu_m,\Pic^{0,\red}_X)=0$ by \cite[Cor.\ 17.5, II 14.2]{oort}.
A diagram chase shows that $f$ is injective, and
the right vertical map is an isomorphism.
The Corollary follows because $\colim \Hom_{\bar k}(\mu_m,A_X)$
is divisible and  $\colim \Hom_{\bar k}(\mu_m,NS_X)$ is finite.
\proofend

The above result should be compared to \cite[Prop.6.2]{ichtor},
where we show that, for every proper scheme over an algebraically closed
field, the higher Chow group of zero-cycles
$CH_0(X,1,\Z/m)$ is the Pontrjagin dual of $H^1_\et(X,\Z/m)$.
This implies a short exact sequence
$$ 0\to {}_\tor A_X^t(k) \to CH_0(X,1,\Q/\Z) \to \chi(T_X)\otimes \Q/\Z\to 0,$$
for $A_X^t$ the dual abelian variety of $A_X$, and $\chi(T_X)$
the character module of $T_X$. However, in this case the contribution
from the torus and from the abelian variety are not compatible with the
coefficient sequence
$$ 0\to CH_0(X,1)\otimes\Q/\Z \to CH_0(X,1,\Q/\Z) \to {}_\tor CH_0(X)\to 0$$
as in Corollary \ref{detail}.

Looking at tangent spaces, the previous Corollary gives a
dimension formula:

\begin{corollary}
Let $l$ be a prime different from $\chr k$. Then
$$ \dim_k H^1(X,\O_X)=
\dim U_X + \dim_k Lie(NS^0_X)+ \rk H^1_\et(X,\Z)+
\textstyle\frac{1}{2} \corank_l H^1_\et(\bar X,\Q_l/\Z_l).$$
\end{corollary}

\section{The unipotent part}
Let $N\Pic(X):= \ker\big(\Pic(X[t])\xrightarrow{0^*} \Pic (X)\big)$.
Since $t\mapsto 0t$ induces $x\mapsto 0x$ on the typical part,
$\Pic(X[t])_{[1]}$ is a subgroup of $N\Pic(X)$. In
\cite{weibelcontracted}, Weibel shows that for every scheme there is
a direct sum decomposition
$$\Pic(X[t,t^{-1}])\cong\Pic(X)\oplus N\Pic(X)\oplus N\Pic(X)\oplus
H^1_\et(X,\Z).$$

\proof (Theorem \ref{npic}).
We show first that $N\Pic(X)=\ker\big(  U_X(\A^1)\to U_X(k)\big)$.
Since there are no non-trivial morphisms of schemes from
$\A^1_k$ to an abelian variety, a torus, an infinitesimal group,
or a discrete group, we
see that the kernel of $U_X(\A^1_k)\to U_X(k)$ agrees
with the kernel of $\Pic_X(\A^1_k)\to \Pic_X(k)$.
Let $p:X\to k$ and $p':X\times \A^1_k\to \A^1_k$ be the structure morphisms.
Then the Leray spectral sequence gives a commutative diagram
$$\begin{CD}
0@>>>H^1_\et(\A_k^1,p_*'\G_m)@>>> \Pic(X\times\A^1_k)@>>> \Pic_X(\A^1_k)
@>>> H^2_\et(\A_k^1,p'_*\G_m)\\
@III @VVV @VVV @VVV @VVV\\
0@>>>H^1_\et(k,p_*\G_m)@>>> \Pic(X)@>>>  \Pic_X(k)@>>> H^2_\et(k,p_*\G_m),
\end{CD}$$
and it suffices to show that the outer vertical maps are isomorphisms.
Let $X\xrightarrow{g} L\to k$ be the Stein factorization of $p$, such that
$\O_L\cong g_*\O_X$ and $L$ is the spectrum of an Artinian $k$-algebra.
Since $\A^1_k\to k$ is flat,
$p'_*\O_{X\times\A^1_k}=\O_{\A^1_k}\otimes_kp_*\O_X$, and
$X\times \A^1_k\xrightarrow{g'} \A^1_L\ \xrightarrow{} \A^1_k$
is the Stein factorization of $p'$. We obtain
$$H^i_\et(\A_k^1,p'_*\G_m)\cong H^i_\et(\A_L^1,g'_*\G_m)\cong
H^i_\et(\A_L^1,\G_m),$$
and $H^i_\et(k,p_*\G_m)\cong H^i_\et(L,\G_m)$. Hence the terms on the left
vanish because $\Pic(L)=\Pic(\A_L^1)=0$.
To show that $H^2_\et(\A_L^1,\G_m)\to H^2_\et(L,\G_m)$
is an isomorphism, we can assume that $L$ is a local Artinian $k$-algebra
with (perfect) residue field $k'$. By \cite[III Rem.3.11]{milnebook}
we are reduced to showing that
$H^2_\et(\A_{k'}^1,\G_m)\to H^2_\et(k',\G_m)$ is an isomorphism,
and this can be found in \cite[IV Ex.2.20]{milnebook}.

Given an element $x$ of $N\Pic(X)$, the condition $x\in \Pic(X[t])_{[1]}$
implies that the corresponding $f\in \Hom_{Sch}(\A^1,U_X)$
satisfies $f(nx)=nf(x)$ for all $n$.
If $k$ has characteristic $0$, then $U_X\cong \G_a^r$
for some $r$, and the map $f: \G_a\to U_X$ corresponds to a morphism
of Hopf algebras $f^*:k[x_1,\cdots ,x_r]\to k[t]$. If
$f^*(x_i)=\sum_ja_jt^j$, then
$$\sum_ja_j(nt)^j= nf^*(x_i)= f^*(nx_i)=n\sum_ja_jt^j$$
only if $n^j=n$ for all $n$, hence $j=1$.
\proofend

\example If $k$ has characteristic $p$, then $t\mapsto t^{2p-1}$
induces a map $\G_a\to \G_a$ which is compatible with multiplication
by $n$, but not a homomorphism of group schemes.

\medskip

\begin{corollary}
We have $U_X=0$ if and only if $N\Pic(X)=0$.
\end{corollary}

\proof
This follows from $N\Pic(X)= \ker \big(U_X(\A^1)\to U_X(k)\big)$, because
any unipotent, connected, smooth affine group
is an affine space as a scheme, hence admits a non-trivial
morphism from $\A^1$ which sends $0$ to $0$ if it is non-trivial.
\proofend

The kernel and cokernel of $\Pic_X\to \Pic_{X^\red}$ has been described in
\cite{oortthesis},
hence we will from now assume that $X$ is reduced.
If $X^\sn$ is the semi-normalization of $X$, then the map
$\O_X\to \O_{X^\sn}$ is an injection of sheaves on the same topological space.
For $X^\sn$ reduced and semi-normal, $N\Pic(X^\sn)=0$ by Traverso's theorem
\cite{traverso} together with \cite[Thm. 4.7]{weibelcontracted}.
Hence the Corollary implies that $U_{X^\sn}=0$, and that
$$ U_X=\ker( \Pic_X^{0,\red}\to \Pic_{X^\sn}^{0,\red}).$$
(For curves, this recovers \cite[Prop.9.2/10]{blr}.)
Indeed, by Corollary \ref{detail}, the map
$\Pic_X^{0,\red}\to \Pic_{X^\sn}^{0,\red}$ induces an isomorphism
on the torus and abelian variety part, because it induces an
isomorphism on etale cohomology.

\proof (Theorem \ref{unipotent})
a) We have isomorphims $H^i_\et(X,\Z)\cong H^i_\et(X^\sn,\Z)$, which combined 
with Corollary \ref{detail} shows that the
canonical map $\Pic_X^{0,\red}\to \Pic_{X^\sn}^{0,\red}$ induces an
isomorphism on the torus components, and is an isogeny with kernel
a unipotent group scheme $P_X$ on the abelian variety parts.
\footnote{O.Gabber \cite{gabber} showed that, conversely, any finite unipotent
commutative group scheme can appear as $P_X$.} 
Hence the map is
surjective and the kernel $K_X$ is an extension of $P_X$ by $U_X$. 
Applying the Proposition to the exact sequence of etale sheaves 
$$ 0\to p_*\G_{m,X} \to p_*\G_{m,X^\sn}\to p_*(\G_{m,X^\sn}/\G_{m,X}) \\
\to \Pic_X\to \Pic_{X^\sn}
$$
on $\Spec k$, we obtain the diagram with exact columns
$$\begin{CD}
0@>>> K_X @>>> \Pic_X^{0,\red} @>>>  \Pic_{X^\sn}^{0,\red}@>>>0 \\
@. @VuVV@VvVV@VVV\\
0@>>> p_*(\G_{m,X^\sn}/\G_{m,X})@>>> \Pic_X @>>>  \Pic_{X^\sn}\\
@. @VVV@VVV@VVV\\
0@>>> \coker u @>>> \NS_X @>>>  \NS_{X^\sn} .
\end{CD}$$
Since $v$ is injective, so is $u$. 
The Neron-Severi group schemes are extensions of finitely generated
\'etale group schemes by a finite connected group scheme. 
The isomorphism $H^2_\et(X,\mu_m)\cong H^2_\et(X',\mu_m)$ implies that
$\Pic_X(\bar k)/m\to \Pic_{X^+}(\bar k)/m$ is injective for any 
$m$ prime to $p$, and since $\Pic_X^{0,\red}(\bar k)$ and 
$\Pic_{X^\sn}^{0,\red}(\bar k)$ are $m$-divisible, the same holds for  
$\NS_X(\bar k)/m\to \NS_{X^+}(\bar k)/m$, and consequently for 
$\NS_X[\frac{1}{p}]\to \NS_{X^+}[\frac{1}{p}]$.
Thus $\coker u$ is contained in the
extension of the $p$-primary torsion subgroup $\NS_X\{p\}$ by 
the finite connected group scheme $\NS_X^0$.
Finally, the isomorphism
$\Hom_{\bar k}(\mu_m,\Pic_X)\cong H^1_\et(\bar X,\Z/m)$
from \cite[III Prop.\ 4.16]{milnebook} together with the isomorphism
$H^i_\et(X,\Z)\cong H^i_\et(X^\sn,\Z)$ and the result on $K_X$ 
shows that the three right maps in the diagram
induce isomorphisms on $\Hom_{\bar k}(\mu_m , - )$ for  
all $m$, hence the three groups on the left are unipotent.

b) Recall that $q:X_T\to T$, and consider the diagram
$$\begin{CD}
0@>>> H^1_\et(T,q_*\G_{m,X\times T})@>>> \Pic(X\times T)@>>> \Pic_{X/k}(T)
@>>> H^2_\et(T,q_*\G_{m,X\times T})\\
@III@VVV @VrVV @VsVV @VVV \\
0@>>>H^1_\et(T,q_*\G_{m,X^\sn\times T})@>>>\Pic(X^\sn\times T)
@>>> \Pic_{X^\sn/k}(T)@>>> H^2_\et(T,q_*\G_{m,X^\sn\times T}).
\end{CD}$$
Since $q_*\G_{a,X\times T}=H^0(X,\O_X)\otimes\O_T=
H^0(X^\sn,\O_{X^\sn})\otimes\O_T=q_*\G_{a,X^\sn\times T}$
is an isomorphism as in Proposition \ref{bef}a),
the outer maps are isomorphisms, and it suffices to calculate $\ker r$.
Let $Y=X\times T$ and $Y'=X^\sn\times T$, and consider
the tautological map
$$f:\{\O_{Y}\text{-line bundles }{\mathcal L}\subseteq \O_{Y'}
\text{ which are invertible in }\O_{Y'}\} \to  \Pic(Y).$$
It suffices to show the following statements:
\begin{enumerate}
\item The image of $f$ is contained in
$\ker\big(\Pic(Y)\to \Pic(Y')\big)$.
\item $f$ surjects onto $\ker\big( \Pic(Y)\to \Pic(Y')\big)$.
\item $f$ is injective.
\end{enumerate}

a) We claim that the map
$\L\otimes_{\O_Y}\O_{Y'}\to \O_{Y'}\otimes_{\O_Y}\O_{Y'}
\xrightarrow{\mu} \O_{Y'}$ is an isomorphism. We can check this
on an affine covering, and in this case it
is proved in \cite[Lemma 2.2(4)]{robertssing}.

b) Let $\L\in \Pic(Y)$ with
$\L\otimes_{\O_Y}\O_{Y'}\cong \O_{Y'}$. Since $\L$ is flat,
we get an injection $\L=\L\otimes_{\O_Y}\O_{Y}\to
\L\otimes_{\O_Y}\O_{Y'}\cong \O_{Y'}$. We claim that the inverse
of $\L$ in $\O_{Y'}$ is the sheaf associated to the presheaf
$U\mapsto \{x\in \O_{Y'}(U) | x\L(U)\subseteq \O_Y(U)\}\subseteq\O_{Y'}(U)$.
This can be checked on an affine covering, and then it is
\cite[Lemma 2.2(2)]{robertssing}.

c) Let $\L$ and $\L'$ be subsheaves of $\O_{Y'}$ which
are invertible in $\O_{Y'}$ and isomorphic as abstract invertible sheaves.
Multiplying with the inverse of $\L'$ inside $\O_{Y'}$, it suffices
to show that if $\L$ is a subsheaf of $\O_{Y'}$, and $f:\O_Y\to \L$
an isomorphism, then $\L=\O_Y\subseteq \O_{Y'}$. But $f(1)$ is a
global unit of $\O_{Y'}(Y)$, and by Proposition \ref{bef}a),
$\O_{Y}(Y)^\times =\O_{Y'}(Y)^\times $. Hence $\L= f(1)^{-1}\L=\O_Y $.
\proofend

\end{document}